\documentclass{article}

\usepackage{microtype}
\usepackage{graphicx}
\usepackage{subfigure}
\usepackage{booktabs} 
\usepackage{caption}

\usepackage{amsthm}
\usepackage{amssymb,amsmath}

\newcommand{\mbfx}{\mathbf{x}}

\newcommand{\mbfy}{\mathbf{y}}
\newcommand{\mbfz}{\mathbf{z}}

\newcommand{\mbfw}{\mathbf{w}}

\newcommand{\mbfu}{\mathbf{u}}
\newcommand{\mbfv}{\mathbf{v}}

\newcommand{\mbfa}{\mathbf{a}}

\newcommand{\mbr}{\mathbb{R}}

\newtheorem{theorem}{Theorem}
\newtheorem{proposition}{Proposition}
\newtheorem{remark}{Remark}

\newcommand{\dist}{\text{dist}}
\DeclareMathOperator*{\argmin}{arg\,\min}

\usepackage{hyperref}

\usepackage[accepted]{icml2018}

\icmltitlerunning{A DCA-Like Algorithm and its Accelerated Version}

\usepackage{multirow}
\usepackage{tabularx}

\begin{document}

\twocolumn[
\icmltitle{A DCA-Like Algorithm and its Accelerated Version with Application in Data Visualization}



\icmlsetsymbol{equal}{*}

\begin{icmlauthorlist}
\icmlauthor{Hoai An Le Thi}{lgipm}
\icmlauthor{Hoai Minh Le}{lgipm}
\icmlauthor{Duy Nhat Phan}{lgipm}
\icmlauthor{Bach Tran}{lgipm}
\end{icmlauthorlist}

\icmlcorrespondingauthor{Duy Nhat Phan}{duy-nhat.phan@univ-lorraine.fr}

\icmlaffiliation{lgipm}{Department Informatics and Application, LGIPM, University of Lorraine, France}

\icmlkeywords{Machine Learning, ICML}

\vskip 0.3in
]




\begin{abstract}
In this paper, we present two variants of DCA (Different of Convex functions Algorithm) to solve the constrained sum of differentiable function and composite functions minimization problem, with the aim of increasing the convergence speed of DCA. In the first variant, DCA-Like, we introduce a new technique to iteratively modify the decomposition of the objective function. This successive decomposition could lead to a better majorization and consequently a better convergence speed than the basic DCA. We then incorporate the Nesterov's acceleration technique into DCA-Like to give rise to the second variant, named Accelerated DCA-Like. The convergence properties and the convergence rate under Kudyka-Lojasiewicz assumption of both variants are rigorously studied. As an application, we investigate our algorithms for the t-distributed stochastic neighbor embedding. Numerical experiments on several benchmark datasets illustrate the efficiency of our algorithms.
\end{abstract}

\section{Introduction}\label{Sec:intro}
\label{introduction}
In this work, we consider the constrained sum of differentiable function and composite functions minimization problem of the form
\begin{equation}\label{model}
\min_{\mbfx\in X}\left\{F(x) = f(\mbfx) + \sum_{i=1}^m h_i(g_i(\mbfx_i))\right\},
\end{equation}
where $f:\mathbb{R}^n \rightarrow\mathbb{R}$ is a continuously differentiable (possibly nonconvex) function with $L$-Lipschitz continuous gradient; $g_i:\mathbb{R}^{n_i} \rightarrow\mathbb{R}$ ($i=1\ldots n$) are continuous convex functions (possibly nonsmooth) with $\sum_{i=1}^mn_i = n$; $h_i$ are concave increasing and $\partial (-h_i)(t) \subset \mathbb{R}_-$ if $t \geq g_i(\mbfx_i)$; $X$ is a closed convex subset of $\mbr^n$.

The assumptions on $f$, $g$ and $h$ are sufficiently large to cover numerous classes of problem arising from several domains such as Machine Learning, computational biology, image processing, etc. For instance, numerous problems in Machine Learning are formulated as a minimization of the trade-off of a loss function $f(x)$ and a regularizer function $r(x)$, i.e. $\min\limits_{x} f(x) + \lambda r(x)$ with trade-off parameter $\lambda>0$. Since we consider $f$ as possibly nonconvex, it covers several loss functions such as least square, squared hinge, logistic loss, and many other loss functions. Let us now show that the sum of composite functions $h_i(g_i(x))$ can cover numerous existing regularizer functions $r(x)$.

$\bullet$ If $r(x)$ is a convex regularizer such as $\Vert.\Vert_1$, $\Vert.\Vert_2$ or $\Vert.\Vert_\infty$ then we can simply choose $h$ as a linear function $h_i(t) \equiv h(t)=\lambda t$ and $g_i(x)=r(x)$.

$\bullet$ Consider now $r(x)$ as the zero norm ($\Vert.\Vert_0$) usually used in variable selection problem. The zero-norm can be then approximated by a nonconvex regularizers such as: capped-$\ell_1$, exponential function, logarithm function, SCAD, $\ell_p (p<0)$ and $\ell_p (0<p<1)$. The readers are referred to \cite{LT15} for an extensive overview of these nonconvex regularizers. For instance, the capped-$\ell_1$ of a vector $x\in\mathbb{R}^n$ is defined by $r_{cap}(x) = \sum\limits_{i=1}^{n}\min\{1,\theta |x_i|\}$. By defining $h_i(t) \equiv  h(t) = \lambda\min\{1,\theta t\}$ and $g_i(x) = |x_i|$, the capped-$\ell_1$ regularizer takes the form of the sum of composites functions. Similarly, it is easy to show that all aforementioned nonconvex regularizer functions can be expressed as a sum of composite function $h_i(g_i(x))$. 

$\bullet$ In the same way, we can prove that the sum of composite functions covers most of convex ($\ell_{1,1}$, $\ell_{1,2}$, \ldots)  as well as nonconvex mixed-norm regularizers (e.g., $\ell_{q,0}$ with $q=1,2$ or $+\infty$) that are usually used in group variable selection.

\noindent\textit{Paper's contribution.} 
In this work, we investigate new methods based on DC (Difference of Convex functions) programming and DCA (DC Algorithm) to solve \eqref{model}. DCA was introduced in 1985 by T. Pham Dinh in the preliminary state, and extensively developed by H.A. Le Thi and T. Pham Dinh since 1994 to become now classic and increasingly popular (\citet{letthe,letabr2018,PLT98,phacon,pharec} and references therein). DCA has been successfully applied to various nonconvex/nonsmooth programs thanks to its versatility, flexibility, robustness, inexpensiveness and their adaptation to the specific structure of considered problems. The contributions of the paper are multiple, from both theoretical and computation aspects.

By exploiting the special structure of the problem \eqref{model} we will prove that it can be equivalently reformulated as a DC program. Hence DCA can be applied to solve \eqref{model}.

We propose a variant of DCA, named DCA-Like, for accelerating its convergence speed. In fact, basic DCA scheme for solving \eqref{model} requires to compute a parameter $\mu$ greater or equal to the $L$-Lipschitz constant of $f$. In practice, it is difficult (or even impossible) to efficiently compute the $L$-Lipschitz constant. Hence one usually estimates $L$ by a quite large value. However, a large value of $L$ could lead to a low convergence speed of DCA. Different DCA with $\mu$ updating procedure have been developed to deal with this issue \cite{HOAIAN2014388,pharec,phan17}. Theses $\mu$ updating procedures consist in fixing an initial value of $\mu$ and iteratively updating it under some conditions. However, the convergence of these DCAs with $\mu$ updating procedure cannot be proved if we can not estimate an upper bound of the constant $L$. In DCA-Like, we assume that the constant $L$ is unknown and propose a new technique to update the parameter $\mu$ that could lead to a better majorization and consequently a better convergence speed. By iteratively modifying the parameter $\mu$, we also modify the decomposition of the objective function. Note that, by keeping $\mu$ as small as possible, we cannot ensure that the successive decompositions of the objective function are DC decomposition. However, we can prove that the 	convergence of DCA-Like is still guaranteed. 

To further speed up the convergence of DCA, we propose a second variant of  DCA, named Accelerated DCA-Like, by incorporating an acceleration technique based on a linear extrapolation into DCA-Like. 

We study the theoretical convergence of the proposed algorithms. DCA-Like and Accelerated DCA enjoy all the convergence properties of DCA. Furthermore, we prove that every limit point of the sequence generated by DCA-Like and Accelerated DCA-Like is a critical point of \eqref{model}. In particular, building on the powerful Kudyka-Lojasiewicz property, we show that each bounded sequence generated by DCA-Like globally converges to a critical point. We also prove their stronger results on convergence rate under the Kudyka-Lojasiewicz assumption. 

Finally, to evaluate the performance of DCA-Like and Accelerated DCA-Like, we consider the t-distributed Stochastic Neighbor Embedding (t-SNE). The t-SNE, a dimensionality reduction
algorithm \cite{Maaten2008} originally introduced for data visualizing, has been widely used in various applications, e.g. computer security, bio-informatics, etc. DCA-Like and ADCA-Like applied to the t-SNE are inexpensive : the solution of the convex sub-problem can be explicitly computed. Furthermore, we prove that, Majorization Minimization \cite{Yang2015}, the best state-of-the-art algorithm for t-SNE is nothing else but DCA-Like applied to the t-SNE model. We carefully conduct the numerical experiments and provide a comparison of proposed algorithms on several benchmark datasets. 

The remainder of the paper is organized as follows. An overview of DC programming and DCA are given in Section \ref{Sec:DCA}. In Section \ref{Sec:DCA-Like}, we introduce DCA-Like and study its convergence properties as well as its convergence rate. Accelerated DCA-Like is presented in Section \ref{Sec:ADCA-Like}. The numerical experiments on the t-SNE problem are reported in Section \ref{Sec:experiment}. Due to the space limit, all theorem's proofs are provided in the supplementary document.

\section{Overview of DC programming and DCA} \label{Sec:DCA}
DC programming and DCA constitute the backbone of smooth/nonsmooth nonconvex programming and global optimization. They address the standard DC program
\[
\alpha = \inf \{F(x) := G(x) - H(x) \,\vert\, x \in \mathbb{R}^n\} \quad (P_{dc}),
\]
where $G, H$ are lower semi-continuous proper convex functions on $\mathbb{R}^n$. Such a function $F$ is called a DC function, and $G - H$ a DC decomposition of $F$ while $G$ and $H$ are the DC components of $F$. Note that any convex constrained DC program can be rewritten in the standard form $(P_{dc})$ by using the indicator function on $C$, defined by $\chi _{C}(x)=0$ if $x\in C$, +$\infty $ otherwise.
\begin{equation*}
\begin{array}{ll}
&\inf \{F(x):=G(x)-H(x):x\in C\mathrm{\ }\} \\
=&\inf \{\chi
_{C}(x)+G(x)-H(x):x\in \mathrm{I\!R}^{n}\}.
\end{array}
\end{equation*}%
The modulus of strong convexity of $\theta$ on $\Omega$, denoted by $\mu(\theta,\Omega)$ or $\mu(\theta)$ if $\Omega = \mathbb{R}^n$, is given by
\[
\mu(\theta,\Omega) = \sup \{\mu \ge 0 : \theta - (\mu /2)\|.\|^2 \text{~is convex on} ~\Omega\}.
\]
One says that $\theta$ is \emph{strongly convex} on $\Omega$ if $\mu(\theta,\Omega) > 0$.

For a convex function $\theta$, the subdifferential of $\theta$ at $x_0 \in \text{dom} \theta := \{x \in \mathbb{R}^n: \theta(x_0) < +\infty\}$, denoted by $\partial \theta(x_0)$, is defined by
\[
\partial\theta (x_0) := \{y \in \mathbb{R}^n : \theta (x) \ge \theta (x_0) + \langle x-x_0,y\rangle, \forall x \in \mathbb{R}^n\}.
\]
The subdifferential $\partial \theta (x_{0})$ generalizes the derivative in the sense that $\theta $ is differentiable at $x_{0}$ if and only if $\partial \theta (x_{0}) \equiv \{\nabla_x\theta (x_{0})\}.$

A point $x^*$ is called a \emph{critical point} of $G-H$, or a generalized Karush-Kuhn-Tucker point (KKT)  of (P$_{dc}$)) if
$\partial H(x^*) \cap \partial G(x^*) \ne \emptyset$.

The main idea of DCA is simple: each iteration $k$ of DCA approximates the concave part $-H$ by its affine majorization (that corresponds to taking $y^k\in \partial H(x^k))$ and computes $x^{k+1}$ by solving the resulting convex problem.
\[\min \{G(x) - \langle x,y^k\rangle : x \in \mathbb{R}^n\} \quad (P_k).\]
The sequence $\{x^k\}$ generated by DCA enjoys the following properties \citep{letthe,phacon}:

(i) The sequence $\{F(x^k)\}$ is decreasing.

(ii) If $F(x^{k+1}) = F(x^k)$, then $x^k$ is a critical point of $(P_{dc})$ and DCA terminates at $k$-th iteration.

(iii) If $\mu (G)+\mu (H)>0$ then the series $\{\Vert x^{k+1}-x^{k}\Vert ^{2}\}$ converges.

(iv) If the optimal value $\alpha$ of $(P_{dc})$ is finite and the infinite sequence $\{x^k\}$ is bounded then every limit point of the sequence $\{x^k\}$  is a critical point of $G-H$.

\section{DCA-Like}\label{Sec:DCA-Like}
We first equivalently reformulate the problem \eqref{model} as follows
\begin{equation}\label{equi_model}
\min_{(\mbfx, \mbfz)}\left\{\varphi(\mbfx,\mbfz):=\chi_\Omega(\mbfx,\mbfz) + f(\mbfx) + \sum_{i=1}^m h_i(z_i)\right\},
\end{equation}
where $\Omega = \{(\mbfx,\mbfz): \mbfx\in X, g_i(\mbfx_i) \leq z_i, i=1,...,m\}$. Denote by $g(\mbfx)$ the vector given by $g(\mbfx) = (g_1(\mbfx_1),...,g_m(\mbfx_m))$. The problems \eqref{model} and \eqref{equi_model} are equivalent in the following sense.

\begin{proposition}\label{proposition}
A point $\mbfx^*\in X$ is a global (resp. local) solution to the problem \eqref{model} if and only if $(\mbfx^*,g(\mbfx^*))$ is a global (resp. local) solution to the problem \eqref{equi_model}.
\end{proposition}
In the remaining of this paper, we consider the problem \eqref{equi_model} instead of \eqref{model}. The objective function $\varphi(\mbfx,\mbfz)$ of \eqref{equi_model} can be rewritten as 
\begin{equation}\label{decompsition}
\varphi(\mbfx,\mbfz) = G_\mu(\mbfx,\mbfz) - H_\mu(\mbfx, \mbfz),
\end{equation}
where $G_\mu(\mbfx,\mbfz) := \frac{\mu}{2}\|\mbfx\|^2 + \chi_\Omega(\mbfx,\mbfz)$ and $H_\mu(\mbfx,\mbfz) := \frac{\mu}{2}\|\mbfx\|^2 - f(\mbfx) - \sum_{i=1}^m h_i(z_i)$
 with $\mu>0$. It is easy to see that $G_\mu(\mbfx,\mbfz)$ is convex since $\Omega$ is a convex set. On the other hand, $f$ is differentiable with $L$-Lipschitz constant, hence $\frac{\mu}{2}\|\mbfx\|^2 - f(\mbfx)$ is convex if $\mu\geq L$. Consequently, $H_\mu(\mbfx,\mbfz)$ is convex and \eqref{equi_model} is a DC program with $\mu\geq L$. In the basic DCA scheme applied to \eqref{equi_model}, one needs to determine the constant $L$ and then choose $\mu \geq L$. In practice, $L$ can not be computed efficiently and is usually estimated by a quite large value. However, a large value of $\mu$ could lead to a low convergence rate of DCA. DCA based algorithms with $\mu$ updating procedure have been investigated in different work \cite{HOAIAN2014388,pharec,phan17} to deal with this issue. For instance, starting with a small value of $\mu$, one increases $\mu$ if the objective value increases in DCA scheme ($\varphi(\mbfx^{k+1},\mbfz^{k+1}) > \varphi(\mbfx^{k},\mbfz^{k})$), i.e. $\mu$ is not large enough to ensure the convexity of $H_\mu(\mbfx,\mbfz)$. One can also start with a large value of $\mu$ and then decreases it as long as the objective function is decreasing. 

In this work, we propose another technique to update $\mu$ based on another criterion. More precisely, at each iteration $k$, we choose $\mu_k$ by finding the smallest number $i_k$ such that with $\mu_k = \eta^{i_k}\max\{\mu_0,\delta\mu_{k-1}\}$ ($\eta>1,0 < \delta < 1$)
\begin{equation}\label{search}
\begin{array}{lll}
H_{\mu_k}(\mbfx^{k+1},\mbfz^{k+1})& \geq &H_{\mu_k}(\mbfx^k,\mbfz^k) + \langle \mbfy^k, \mbfx^{k+1}-\mbfx^k\rangle\\& + &\langle \xi^k, \mbfz^{k+1}-\mbfz^k\rangle,
\end{array}
\end{equation}
where $(\mbfy^k,\xi^k)\in\partial H_{\mu_k}(\mbfx^k,\mbfz^k)$ and $(\mbfx^{k+1},\mbfz^{k+1})$ minimizes the following convex problem
\begin{equation}\label{convex1}
\min\left\{G_{\mu_k}(\mbfx,\mbfz) - \langle \mbfy^k, \mbfx\rangle - \langle \xi^k, \mbfz\rangle\right\}.
\end{equation}
The convex sub-problem \eqref{convex1} can be rewritten as follows
\begin{equation}\label{convex2}
\min_{(\mbfx,\mbfz)\in \Omega}\left\{\frac{\mu_k}{2}\|\mbfx\|^2 - \langle \mbfy^k, \mbfx\rangle + \sum_{i=1}^m (- \xi^k_i) \mbfz_i\right\},
\end{equation}
where $\xi^k_i\in\partial (- h_i)(z^k_i)$. Since $-\xi^k_i \geq 0$, the solution $(\mbfx^{k+1},\mbfz^{k+1})$ to the problem \eqref{convex2} is given by
\begin{equation}
\begin{aligned}
\mbfx^{k+1}&=\argmin_{\mbfx\in X}\{\frac{\mu_k}{2}\|\mbfx\|^2 - \langle \mbfy^k, \mbfx\rangle + \sum_{i=1}^m (- \xi^k_i) g_i(\mbfx_i)\},\\
z^{k+1}_i &= g_i(\mbfx^{k+1}_i), \ i = 1,...,m.
\end{aligned}
\end{equation}
DCA-Like for solving \eqref{equi_model} is described in Algorithm \ref{Alg:DCA-Like}.

\begin{algorithm}[tbh]
\caption{DCA-Like for solving \eqref{equi_model}}
\begin{algorithmic}\label{Alg:DCA-Like}
   \STATE {\bfseries Initialization:} Choose $\mbfx^0$, $\eta>1, 0<\delta<1$, $\mu_0>0$ and $k\leftarrow 0$. 
   \REPEAT 
   	\STATE 1: Compute $\xi^k_i\in\partial (-h_i)\left(g_i(\mbfx_i^k)\right)$ and $\nabla f(\mbfx^k)$.
    \STATE 2: Set $\mu_k = \max\{\mu_0,\delta\mu_{k-1}\}$ if $k>0$.
    \STATE 3: Compute $\mbfx^{k+1}$ by 
    \begin{equation}\label{DCA-Like-SubPrb}
    \min_{\mbfx\in X}\{\frac{\mu_k}{2}\|\mbfx-\mbfx^k\|^2 + \langle \nabla f(\mbfx^k), \mbfx\rangle + \sum_{i=1}^m (- \xi^k_i) g_i(\mbfx_i)\}.
    \end{equation}
   
   \STATE 4: \textbf{While} $H_{\mu_k}(\mbfx^{k+1},g(\mbfx^{k+1})) < H_{\mu_k}(\mbfx^k,g(\mbfx^k)) + \langle \mu_k \mbfx^k - \nabla f(\mbfx^k), \mbfx^{k+1}-\mbfx^k\rangle + \langle \xi^k, g(\mbfx^{k+1})-g(\mbfx^k)\rangle$ \textbf{do}
   \STATE\quad$\bullet$ $\mu_k \leftarrow \eta\mu_k$.
   \STATE\quad$\bullet$ Update $\mbfx^{k+1}$ by STEP 3.
   \STATE\quad \textbf{End While}
   
    \STATE 5: $k\leftarrow k+ 1$.
    \UNTIL{Stopping criterion.}
\end{algorithmic}
\end{algorithm}

\begin{remark}\label{remark1}
$\bullet$ It is easy to show that the while loop in STEP 4 stops after finitely steps. Indeed, it follows from the convexity of $-h_i$ that for $i=1,...,m$
\begin{equation}\label{2:alg3}
- h_i(g_i(\mbfx^{k+1}))\geq -h_i(g_i(\mbfx^k)) + \langle \xi^k_i,g_i(\mbfx^{k+1}) -g_i(\mbfx^k)\rangle.
\end{equation}
Since $f$ is $L$ - Lipschitz gradient, for $\mu_k\geq L$, we have
\begin{equation}\label{3:alg3}
\begin{array}{ll}
f(\mbfx^k) &+ \langle \nabla f(\mbfx^k),\mbfx^{k+1} - \mbfx^k\rangle\\
 &+ \frac{\mu_k}{2}\|\mbfx^{k+1} -\mbfx^k\|^2 \geq f(\mbfx^{k+1}).
\end{array}
\end{equation}
Summing inequalities \eqref{2:alg3} and \eqref{3:alg3} implies that the inequality \eqref{search} holds. From this, there also exists $\beta>0$ such that $\mu_k\leq \beta L$ for all $k$.

$\bullet$ The backtracking condition \eqref{search} does not imply that $\mu$ is large enough to ensure the convexity of $H_{\mu_k}$. However, we will prove that the convergence properties of DCA-Like are still guaranteed. Moreover, by keeping $\mu$ as small as possible, we can get a closer majorization of $\varphi$, which could lead to a faster converge and better solution.

$\bullet$ We have equivalently reformulated the problem \eqref{model} as a constrained problem \eqref{equi_model} by adding variables $z_i$. According to Algorithm \ref{Alg:DCA-Like}, DCA-Like for \eqref{equi_model} consists in solving the sequence of convex problems \eqref{DCA-Like-SubPrb}. As we can see, the sub-problem only involves the variable $x$.
\end{remark}

\subsection{Convergence analysis of DCA-Like}
In this subsection, we study the convergence of DCA-Like. 
Our first result provides the behavior of the limit points of the sequence $\{x^k\}$ generated by DCA-Like.
\begin{theorem}\label{theorem1}
Let $\{\mbfx^k\}$  be the sequence generated by Algorithm \ref{Alg:DCA-Like}. The following statements hold.

(i) The sequence $\{\varphi(\mbfx^k,g(\mbfx^k))\}$ is decreasing. More precisely, we have
\begin{equation*}\label{descent}
\varphi(\mbfx^k,g(\mbfx^k)) - \varphi(\mbfx^{k+1},g(\mbfx^{k+1})) \geq \frac{\mu_k}{2}\|\mbfx^{k+1}-\mbfx^{k}\|^2.
\end{equation*}
(ii) If $\alpha = \inf \varphi(\mbfx,\mbfz) > -\infty$ then
$\sum_{k=0}^{+\infty}\|\mbfx^{k+1}-\mbfx^k\|^2< +\infty,$
and therefore $\lim_{k\rightarrow +\infty}\|\mbfx^{k+1}-\mbfx^k\| = 0$.

(iii) If $\alpha = \inf \varphi(\mbfx,\mbfz) > -\infty$, then any limit point of $\{(\mbfx^k,g(\mbfx^k))\}$ is a critical point of \eqref{equi_model}.
\end{theorem}
Next we study the convergence of the sequence generated by DCA-Like under Kurdyka-Lojasiewicz (KL) assumption. Let $\eta\in(0,\infty]$. Denote by $\mathcal{M}_\eta$ the class of continuous concave functions $\psi:[0,\eta) \rightarrow [0,\infty)$ verifying 

(i) $\psi(0) = 0$ and $\psi$ is continuously differentiable on $(0,\eta)$,

(ii) $\psi'(t)>0$ for all $t\in(0,\eta)$.

Recall that a lower semicontinuous function $\sigma$ satisfies the KL property \citep{attpro} at $\mbfu^*\in \text{dom}\ \partial^L \sigma$ if there exists $\eta >0$, a neighborhood $\mathcal{V}$ of $\mbfu^*$ and $\psi\in \mathcal{M}_\eta$ such that for all $u\in\mathcal{V}\cap\{u:\sigma(\mbfu^*)<\sigma(\mbfu)<\sigma(\mbfu^*)+\eta\}$, one has
\begin{equation*}
\psi'(\sigma(\mbfu) - \sigma(\mbfu^*))\dist(0,\partial^L\sigma(\mbfu)) \geq 1.
\end{equation*}
Here $\partial^L\sigma(\mbfu)$ denotes the limiting-subdifferential of $\sigma$ at $\mbfu$ \citep{Mor06}. The class of functions $\sigma$ verifying the KL property at all points in dom $\partial^L \sigma$ is very ample, for example, semi-algebraic, subanalytic, and log-exp functions. In particular, these classes of functions satisfy the KL property with $\psi(s) = cs^{1-\theta}$, for some $\theta\in[0,1)$ and $c>0$.

In the theorem below, we provide sufficient conditions that guarantee the convergence of the whole sequence $\{\mbfx^k\}$ generated by DCA-Like. These conditions include  the KL property of $\varphi$ and the differentiability with locally Lipschitz derivative of $h_i$. Moreover, if the function $\psi$ appearing in the KL inequality has the form $\psi(s) = cs^{1-\theta}$ with $\theta\in[0,1)$ and $c>0$, then we obtain the rates of convergence for the both sequences $\{\mbfx^k\}$ and $\{\varphi(\mbfx^k,g(\mbfx^k))\}$.

\begin{theorem}\label{theorem2}
Suppose that $\inf\varphi(\mbfx,\mbfz)>-\infty$ and $h_i$ is differentiable with locally Lipschitz derivative. Assume further that $\varphi$ has the KL property at any point $(\mbfx,\mbfz)\in \text{dom}\ \partial^L \varphi$. If $\{\mbfx^k\}$ generated by DCA-Like is bounded, then the whole sequence $\{\mbfx^k\}$ converges to $\mbfx^*$, which $(\mbfx^*,g(\mbfx^*))$ is a critical point of \eqref{equi_model}. Moreover, if the function $\psi$ appearing in the KL inequality has the form $\psi(s) = cs^{1-\theta}$ with $\theta\in[0,1)$ and $c>0$, then the following statements hold

(i) If $\theta = 0$, then the sequences $\{\mbfx^k\}$ and $\{\varphi(\mbfx^k,g(\mbfx^k))\}$ converge in a finite number of steps to $\mbfx^*$ and $\varphi^*$, respectively.

(ii) If $\theta \in (0,1/2]$, then the sequences $\{\mbfx^k\}$ and $\{\varphi(\mbfx^k,g(\mbfx^k))\}$ converge linearly to $\mbfx^*$ and $\varphi^*$, respectively.

(iii) If $\theta \in (1/2,1)$, then there exist positive constants $\delta_1$, $\delta_2$ and $N_0$ such that
$\|\mbfx^k - \mbfx^*\| \leq \delta_1k^{-\frac{1-\theta}{2\theta-1}}$ and 
$\varphi(\mbfx^k,g(\mbfx^k)) - \varphi^* \leq \delta_2 k^{-\frac{1}{2\theta-1}}$
for all $k \geq N_0$. 
\end{theorem}
 As shown in Theorem \ref{theorem2}, the whole bounded sequence $\{\mbfx^k\}$ converges to $\mbfx^*$. In particular, the both sequences $\{\mbfx^k\}$ and $\{\varphi(\mbfx^k,g(\mbfx^k))\}$ converge in finite iterations when $\theta=0$, converge with a linear rate when $\theta \in (0,1/2]$ and a sub-linear rate when $\theta \in (1/2,1)$.
 
\section{Accelerated DCA-Like}\label{Sec:ADCA-Like}

We now introduce the Accelerated DCA-Like (ADCA-Like) for solving the problem \eqref{equi_model}. According to the DCA-Like scheme, at each iteration, one computes $\mbfx^{k+1}$ from $\mbfx^k$ by solving the convex sub-problem \eqref{DCA-Like-SubPrb}. The idea of ADCA-Like, in order to accelerate the convergence of DCA-Like, is to find a point $\mbfw^k$ which is better than $\mbfx^k$ for the computation of $\mbfx^{k+1}$. 
In this work, we consider $\mbfw^k$ as an extrapolated point of the current iterate $\mbfx^k$ and the previous iterate $\mbfx^{k-1}$:
\begin{equation*}
\mbfw^k = \mbfx^k +  \frac{t_k-1}{t_{k+1}}\left(\mbfx^k-\mbfx^{k-1}\right), 
\end{equation*}
where $t_{k+1} = \frac{1+\sqrt{1+4t^2_{k}}}{2}$. If $\mbfw^k$ is better than the last iterate $\mbfx^k$, i.e., $\varphi(\mbfw^k,g(\mbfw^k)) \leq \varphi(\mbfx^k,g(\mbfx^k))$ then $\mbfw^k$ will be used instead of $\mbfx^k$ to compute $\mbfx^{k+1}$. Note that, ADCA-Like does not require any particular property of the sequence $t$. We choose the above sequence as it has interesting convergence rate \cite{becafa}. The proposed algorithm is described in Algorithm \ref{ADCA-like}.

\begin{algorithm}[tbh]
\caption{ADCA-Like for solving \eqref{equi_model}}
\begin{algorithmic}\label{ADCA-like}
   \STATE {\bfseries Initialization:} Choose $\mbfx^0$, $\mbfw^0=\mbfx^0$, $\eta>1, 0<\delta<1$, $\mu_0>0$ and $k\leftarrow 0$. 
   \REPEAT 
     \STATE 1: \textbf{If} $\varphi(\mbfw^k,g(\mbfw^k)) \leq \varphi(\mbfx^k,g(\mbfx^k))$ 
     \textbf{then} set $\mbfv^k = \mbfw^k$ \textbf{else} set $\mbfv^k = \mbfx^k$. 
   	\STATE 2: Compute $\xi^k_i\in\partial (-h_i)\left(g_i(\mbfv_i^k)\right)$ and $\nabla f(\mbfv^k)$.
    \STATE 3: Set $\mu_k = \max\{\mu_0,\delta\mu_{k-1}\}$ if $k>0$.
    \STATE 4: Compute $\mbfx^{k+1}$ by 
    \begin{equation*}
    \min_{\mbfx\in X}\{\frac{\mu_k}{2}\|\mbfx-\mbfv^k\|^2 + \langle \nabla f(\mbfv^k), \mbfx\rangle + \sum_{i=1}^m (- \xi^k_i) g_i(\mbfx_i)\}.
    \end{equation*}
   
   \STATE 5: \textbf{While} $H_{\mu_k}(\mbfx^{k+1},g(\mbfx^{k+1})) < H_{\mu_k}(\mbfv^k,g(\mbfv^k)) + \langle \mu_k \mbfv^k - \nabla f(\mbfv^k), \mbfx^{k+1}-\mbfv^k\rangle + \langle \xi^k, g(\mbfx^{k+1})-g(\mbfv^k)\rangle$ \textbf{do}
   \STATE\quad$\bullet$ $\mu_k \leftarrow \eta\mu_k$.
   \STATE\quad$\bullet$ Update $\mbfx^{k+1}$ by STEP 4.
   \STATE\quad \textbf{End While}
   \STATE 6: Compute $t_{k+1} = \frac{1+\sqrt{1+4t^2_{k}}}{2}$.
   \STATE 7: Compute $\mbfw^{k+1} = \mbfx^{k+1} + \frac{t_k-1}{t_{k+1}}\left(\mbfx^{k+1}-\mbfx^k\right)$.
   
    \STATE 8: $k\leftarrow k+ 1$.
    \UNTIL{Stopping criterion.}
\end{algorithmic}
\end{algorithm}

\subsection{Convergence analysis of accelerated DCA-Like}
The following theorem shows that any limit point of the sequence generated by accelerated DCA-Like is a critical point of \eqref{equi_model}.
\begin{theorem}\label{theorem3}
Let $\{x^k\}$ be the sequence generated by Algorithm \ref{ADCA-like}. The following statements hold

(i) The sequence $\{\varphi(\mbfx^k,g(\mbfx^k))\}$ is decreasing. More precisely, we have
\begin{equation*}\label{descent}
\varphi(\mbfx^k,g(\mbfx^k)) - \varphi(\mbfx^{k+1},g(\mbfx^{k+1})) \geq \frac{\mu_k}{2}\|\mbfx^{k+1}-\mbfv^{k}\|^2.
\end{equation*}
(ii) If $\alpha = \inf \varphi(\mbfx,\mbfz) > -\infty$ then
$\sum_{k=0}^{+\infty}\|\mbfx^{k+1}-\mbfv^k\|^2< +\infty$
and therefore $\lim_{k\rightarrow +\infty}\|\mbfx^{k+1}-\mbfv^k\| = 0$.

(iii) If $\alpha = \inf \varphi(\mbfx,\mbfz) > -\infty$, then any limit point of $\{(\mbfx^k,g(\mbfx^k))\}$ is a critical point of \eqref{equi_model}.
\end{theorem}
The sufficient descent property (i) of Theorem \ref{theorem3} is different from Theorem \ref{theorem1} due to the intermediate variable $\mbfv^k$. Hence, neither the convergence of the whole sequence $\{\mbfx^k\}$ nor convergence rate for $\{\|\mbfx^k-\mbfx^*\|\}$ can be achieved. However, we can still obtain some exciting results for the sequence $\{\varphi(\mbfx^k,g(\mbfx^k))\}$ under the KL assumption. We have the following theorem.
\begin{theorem}\label{theorem4}
Suppose that $\inf\varphi(\mbfx,\mbfz)>-\infty$ and $h_i$ is differentiable with locally Lipschitz derivative. Assume further that $\varphi$ has the KL property at any point $(\mbfx,\mbfz)\in \text{dom}\ \partial^L \varphi$ with $\psi(s) = cs^{1-\theta}$ for some $\theta\in[0,1)$ and $c>0$. If $\{\mbfx^k\}$ generated by accelerated DCA-Like is bounded, then the following statements hold.

(i) If $\theta = 0$, then the sequence $\{\varphi(\mbfx^k,g(\mbfx^k))\}$ converges in a finite number of steps to $\varphi^*$.

(ii) If $\theta \in (0,1/2]$, then the sequence $\{\varphi(\mbfx^k,g(\mbfx^k))\}$ converges linearly to $\varphi^*$.

(iii) If $\theta \in (1/2,1)$, then there exist positive constants $\delta$ and $N_0$ such that
$\varphi(\mbfx^k,g(\mbfx^k)) - \varphi^* \leq \delta k^{-\frac{1}{2\theta-1}}$ 
for all $k \geq N_0$. 
\end{theorem}

\section{Application to t-SNE in visualizing data}\label{Sec:application}

t-SNE was first introduced by \cite{Maaten2008} as a visualization technique for high dimensional data. 
The obstacle of this approach is due to the nature of high-dimensional space where only small pairwise distances are reliable, thus most techniques only try to model such small pairwise distances in the low embedding space. 
t-SNE is a gaining popular method from the family of stochastic neighbor embedding (SNE) methods \cite{Hinton2003}, operates by retaining local pairwise distances. 
It has been applied in many applications such as bioinformatic \cite{Wilson2015}, cancer research, visualize features in neural networks \cite{Mnih2015}, etc.


The t-SNE problem can be described as the minimization of the divergence between two distributions: (1) a distribution that measures pairwise similarities of the input objects and (2) a distribution that measures pairwise similarities of the corresponding low-dimensional points in the embedding. 
Assume we are given a data set of (high-dimensional) input
objects $\mathcal{D} = \{\mbfa_1,...,\mbfa_n\}$ with $\mbfa_i\in\mathbb{R}^d$. Our aim is to learn a low-dimensional embedding in which each object is represented by a point, $\mathcal{E} = \{\mbfx_1,...,\mbfx_n\}$ with $\mbfx_i\in\mathbb{R}^s$. To this end, t-SNE defines joint probabilities $p_{ij}$ that measure the pairwise similarity between objects $\mbfx_i$ and $\mbfx_j$ by symmetrizing two conditional probabilities as $p_{ij} = \frac{p_{j|i} + p_{i|j}}{2n}$, 
where $p_{j|i} = \frac{\exp(-\|\mbfa_i-\mbfa_j\|^2/2\sigma_i^2)}{\sum_{k\neq i}\exp(-\|\mbfa_i-\mbfa_k\|^2/2\sigma_i^2)}\ $ if $i \neq j$, and $0$ otherwise.

In the embedding subspace $\mathcal{E}$, the similarities between two points $\mbfx_i$ and $\mbfx_j$ are measured using a normalized heavy-tailed kernel. Specifically, the embedding similarity $q_{ij}$ between the two points $\mbfx_i$ and $\mbfx_j$ is computed as a normalized Student-t kernel with a single degree of freedom: $q_{ij} = \frac{(1+\|\mbfx_i-\mbfx_j\|^2)^{-1}}{\sum_{k\neq l}(1+\|\mbfx_k-\mbfx_l\|^2)^{-1}}\  \text{if}\ i \neq j, \text{ and } 0 \text{ otherwise}$.
The locations of the embedding points $\mbfx_i$ are determined by minimizing the Kullback-Leibler divergence between the joint distributions $P$ and $Q$:
\begin{equation}\label{tsne}
\min_{\mbfx}\{F(\mbfx) = KL(P||Q) = \sum_{i\neq j}p_{ij}\log \frac{p_{ij}}{q_{ij}} \}.
\end{equation}
The nonconvex optimization problem \eqref{tsne} has been studied in several works \cite{Maaten2008, Yang2009, Vladymyrov2012}, but the most noticeable was presented in \cite{Yang2015}. \citet{Yang2015} presented and compared Majorization Minimization algorithm (MM) with five state-of-the-arts methods, such as gradient descent, gradient descent with momentum \cite{Maaten2008}, spectral direction \cite{Vladymyrov2012}, FPHSSNE \cite{Yang2009} and Limited-memory Broyden-Fletcher-Goldfarb-Shanno (L-BFGS) \cite{Nocedal1980}. The numerical results showed that MM by outperforms all five state-of-the-art optimization methods.

The objective function $F$ of \eqref{tsne} can be rewritten as follows
\begin{equation*}
\begin{aligned}
F(\mbfx)= \sum_{i\neq j}p_{ij}\log p_{ij} &+ \log(\sum_{i\neq j}(1+\|\mbfx_i-\mbfx_j\|^2)^{-1})\\
& + \sum_{i,j}p_{ij}\log(1+\|\mbfx_i-\mbfx_j\|^2). 
\end{aligned}
\end{equation*}
Let $f(\mbfx) = \sum_{i\neq j}p_{ij}\log p_{ij} + \log(\sum_{i\neq j}(1+\|\mbfx_i-\mbfx_j\|^2)^{-1})$, 
$h_{ij}(t) = p_{ij}\log(1+t)$ and $g_{ij}(\mbfx_i,\mbfx_j) = \|\mbfx_i-\mbfx_j\|^2$. It is obvious that $g_{ij}$ are convex functions, and $h_{ij}$ are concave increasing functions whose derivatives are non-negatives and Lipschitz continuous on $[0,+\infty)$. Moreover, the function $f$ is differentiable with $L$-Lipschitz continuous gradient by the following proposition.
\begin{proposition}\label{lipschitz}
The function $f(\mbfx) = \sum_{i\neq j}p_{ij}\log p_{ij} + \log(\sum_{i\neq j}(1+\|\mbfx_i-\mbfx_j\|^2)^{-1})$ is smooth with Lipschitz gradient, where we can choose a Lipschitz constant $L = 6n\sqrt{s}$.
\end{proposition}
Therefore, the nonconvex problem \eqref{tsne} takes the form of \eqref{model}. Thus, we can investigate DCA-Like and ADCA-Like to solve the problem \eqref{tsne}. Note that both DCA-Like and ADCA-Like are also applicable for other variants of SNE such as SNE \cite{Hinton2003}, Symmetric SNE \cite{Maaten2008}, etc.

According to DCA-Like, from $\mbfx^k$, we have to compute $\xi^k_{ij} = \nabla (-h_{ij})(g_{ij}(\mbfx_i^k,\mbfx_j^k)) = -\frac{p_{ij}}{1+\|\mbfx_i^k-\mbfx_j^k\|^2}$ and $\nabla f(\mbfx^k)$ by
\begin{equation}\label{gradient}
\nabla_{\mbfx_i} f(\mbfx^k) = \sum_{j=1}^n\frac{-4(\mbfx_i^k-\mbfx_j^k)(1+\|\mbfx_i^k-\mbfx_j^k\|^2)^{-2}}{\sum_{l\neq m}(1+\|\mbfx_l^k-\mbfx_m^k\|^2)^{-1}},
\end{equation}
and solve the following convex problem
\begin{equation*}
\min_{\mbfx}\{\frac{\mu_k}{2}\|\mbfx - \mbfx^k\|^2 +  \langle \nabla f(\mbfx^k), \mbfx \rangle + \sum_{i,j} -\xi^k_{ij} \|\mbfx_i-\mbfx_j\|^2\}
\end{equation*}
The solution $\mbfx^{k+1}$ to this problem is given by
\begin{equation}\label{solution}
\mbfx^{k+1} = (2\mathcal{L}_{-\xi^k - (\xi^k)^T}+\mu_kI)^{-1}(-\nabla f(\mbfx^k) +\mu_k\mbfx^k),
\end{equation}
where the matrix $\xi^k$ is defined by the elements $\xi_{ij}^k$ and $\mathcal{L}_A$ denotes the matrix with $(\mathcal{L}_A)_{ij} = -A_{ij}$ if $i\neq j$ and $- A_{ii} + \sum_{l=1}^n A_{il}$ otherwise. We observe that the while loop in Algorithm \ref{Alg:DCA-Like} stops if the following inequality holds
\begin{equation}\label{searchmm}
U_{\mu_k}(\mbfx^{k+1},\mbfx^k) \geq F(\mbfx^{k+1}),
\end{equation}
where $U_{\mu_k}(\mbfx^{k+1},\mbfx^k) = F(\mbfx^k) + \langle \nabla f(\mbfx^k),\mbfx^{k+1} - \mbfx^k\rangle + \frac{\mu_k}{2}\|\mbfx^{k+1} -\mbfx^k\|^2 - \langle \xi^k,g(\mbfx^{k+1}) -g(\mbfx^k)\rangle$. From the update rule \eqref{solution} for $\mbfx^{k+1}$ and this stopping criterion for searching $\mu_k$, we can conclude that MM \citep{Yang2015} for \eqref{tsne} is special version of DCA-Like. In summary, DCA-Like for solving t-SNE problem \eqref{tsne} is described in Algorithm \ref{DCA-like-tsne}.

\begin{algorithm}[tbh]
\caption{DCA-Like for \eqref{tsne}}
\begin{algorithmic}\label{DCA-like-tsne}
   \STATE {\bfseries Initialization:} Choose $\mbfx^0$, $\eta>1, 0<\delta<1$, $\mu_0>0$ and $k\leftarrow 0$. 
   \REPEAT 
   	\STATE 1: Compute $\xi^k_{ij}=-\frac{p_{ij}}{1+\|\mbfx_i^k-\mbfx_j^k\|^2}$ and $\nabla f(\mbfx^k)$ by \eqref{gradient}.
    \STATE 2: Set $\mu_k = \max\{\mu_0,\delta\mu_{k-1}\}$ if $k>0$.
    \STATE 3: Compute $\mbfx^{k+1}$ by \eqref{solution}.
       
   \STATE 4: \textbf{While} $U_{\mu_k}(\mbfx^{k+1},\mbfx^k) < F(\mbfx^{k+1})$ \textbf{do}
   \STATE\quad$\bullet$ $\mu_k \leftarrow \eta\mu_k$.
   \STATE\quad$\bullet$ Update $\mbfx^{k+1}$ by STEP 3.
   \STATE\quad \textbf{End While}
   
    \STATE 5: $k\leftarrow k+ 1$.
    \UNTIL{Stopping criterion.}
\end{algorithmic}
\label{alog.DCA-Like-tSNE}
\end{algorithm}

The ADCA-Like for solving \eqref{tsne} is obtained by adding STEP 1,6 and 7 of Algorithm \ref{ADCA-like} to Algorithm \ref{DCA-like-tsne}.

We recall that all semi-algebraic functions and subanalysis functions satisfy the KL property \cite{attpro}, for examples, real polynomial functions, logarithm function, $\ell_p$-norm with $p\geq 0$. In addition, finite sums, products, generalized inverse, compositions of semi-algebraic functions are also semi-algebraic. This implies that the objective function of \eqref{tsne} satisfies the KL property. Hence, DCA-Like and ADCA-Like for solving \eqref{tsne} enjoy all convergence properties provided in Theorems \ref{theorem1}-\ref{theorem4}.

\section{Numerical experiment}\label{Sec:experiment}

To evaluate the performances of our methods, we perform numerical experiments on six real datasets taken from UCI data repository (\textit{letters}, \textit{shuttle}, \textit{sensorless}, \textit{mnist}, \textit{miniboone} and \textit{covertype}). The comparison are realized on three criteria: the objective value $F(\mbfx)$, the number of iterations and the computation time (measured in seconds). Each experiment is repeated $10$ times, then the final result is the average value of each criterion.

As mentioned before, the t-SNE can also be solved by DCA with the DC decomposition \eqref{decompsition}. For DCA scheme, we have to estimate the $L$-Lipschitz constant of $f$. According to Proposition \ref{lipschitz}, we can choose $L = 6n\sqrt{s}$. This value is clearly too large. Hence, we will incorporate a $\mu$ updating procedure into DCA. We start with a small value of $\mu$ and increase $\mu$ if the objective value increases in DCA scheme. For all algorithms, the initial value of $\mu_0$ is set to be $10^{-6}$.

We follow the same process as described in \cite{Yang2015}. For all datasets, k-Nearest Neighbor (with $k = 10$) is employed to construct $\bar{p}_{ij}$, where $\bar{p}_{ij} = 1$ if data point $j$ (reps. i) is one of k nearest neighbors of data point j (reps. i), and $\bar{p}_{ij} = 0$ otherwise. 
$p_{ij}$ is then computed by $p_{ij} = \frac{\bar{p}_{ij}}{\sum_{k, l}\bar{p}_{kl}}$. 
$\mbfx^0$ is drawn from normal distribution $\mathcal{N}(0, 10^{-8})$ for all methods. Early exaggeration technique \cite{Maaten2008} is deployed for first $20$ iterations with the constant value of $4$. 

For large datasets, Barnes-Hut tree approximation is used for reducing computing cost \cite{Maaten2014}. 
This technique is well-known in Neighbor Embedding problems, which provides a good trade-off small loss in gradients and cost function against huge reduction in computation time. 
We set the parameter $\theta_{\text{Barnes-Hut}} = 0.5$.

Stopping conditions of all algorithms are the same, by either (1) number of iterations exceeds $10000$ or (2) $\| \mbfx^{k} - \mbfx^{k-1} \| / \| \mbfx^{k-1}\| \leq 10^{-8}$. Throughout our experiment, the number of embedding dimension is set to $s = 2$. All experiments are performed on a PC Intel (R) Xeon (R) E5-2630 v4 @2.20 GHz of 32GB RAM.

\begin{table}[tbh]
\caption{Comparative results on datasets. 
Bold values correspond to best results for each dataset, $n$ and $d$ are the number of instances and dimensions respectively. Unit of time is second.}
\label{tbl.exp_result}
\vskip 0.1in
\setlength\tabcolsep{4.75pt}
\begin{tabularx}{\columnwidth}{l|l|rrr}
Dataset & Algorithm & Obj. & Iteration & Time \\ \hline
\textit{letters} & DCA & 1.58 & 1185.67 & 652.80 \\
$n$ = 20000 & DCA-Like & \textbf{1.48} & 164.00 & 149.00 \\
$d$ = 16 & ADCA-Like & \textbf{1.48} & \textbf{89.67} & \textbf{96.80} \\ \hline
\textit{shuttle} & DCA & 1.60 & 3519.67 & 6056.48 \\
$n$ = 58000 & DCA-Like & 1.48 & 304.00 & 899.60 \\
$d$ = 9 & ADCA-Like & \textbf{1.43} & \textbf{143.17} & \textbf{517.30} \\ \hline
\textit{sensorless} & DCA & \textbf{3.18} & 1787.67 & 3232.08 \\
$n$ = 58509 & DCA-Like & 3.21 & 313.33 & 953.75 \\
$d$ = 48 & ADCA-Like & 3.19 & \textbf{142.33} & \textbf{492.06} \\ \hline
\textit{mnist} & DCA & 3.44 & 3893.67 & 13752.75 \\
$n$ = 70000 & DCA-Like & 3.45 & 343.50 & 1481.66 \\
$d$ = 784 & ADCA-Like & \textbf{3.43} & \textbf{187.33} & \textbf{869.02} \\ \hline
\textit{miniboone} & DCA & 3.55 & 3401.00 & 20473.33 \\
$n$ = 130064 & DCA-Like & \textbf{3.53} & 469.00 & 4841.01 \\
$d$ = 50 & ADCA-Like & \textbf{3.53} & \textbf{175.67} & \textbf{1820.06} \\ \hline
\textit{covertype} & DCA & 2.12 & 4013.67 & 68302.53 \\
$n$ = 581012 & DCA-Like & 2.09 & 1223.67 & 35719.52 \\
$d$ = 54 & ADCA-Like & \textbf{1.92} & \textbf{227.67} & \textbf{6875.20}
\end{tabularx}
\end{table}
\begin{table*}[tbh]
\centering
\begin{tabular}{ccc}
(a) \textit{letters} & (b) \textit{shuttle}  & (c) \textit{sensorless} 
\\
\includegraphics[width=0.3\linewidth]{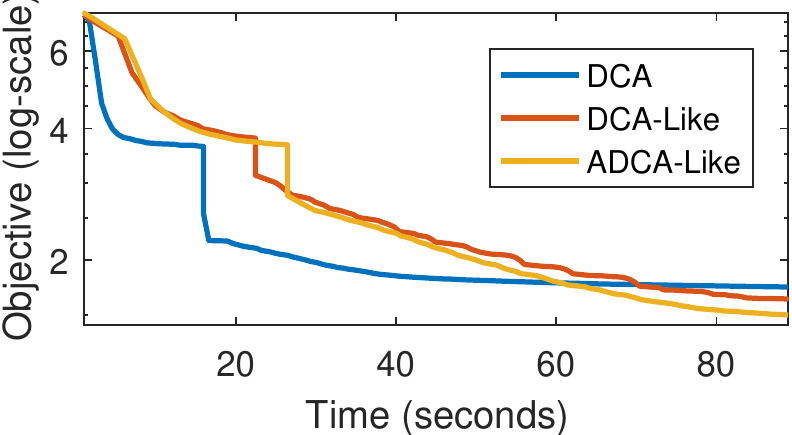} &
\includegraphics[width=0.3\linewidth]{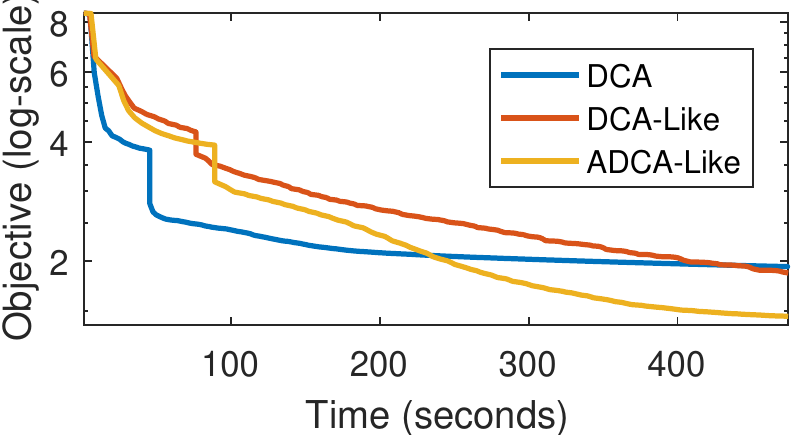} &
\includegraphics[width=0.3\linewidth]{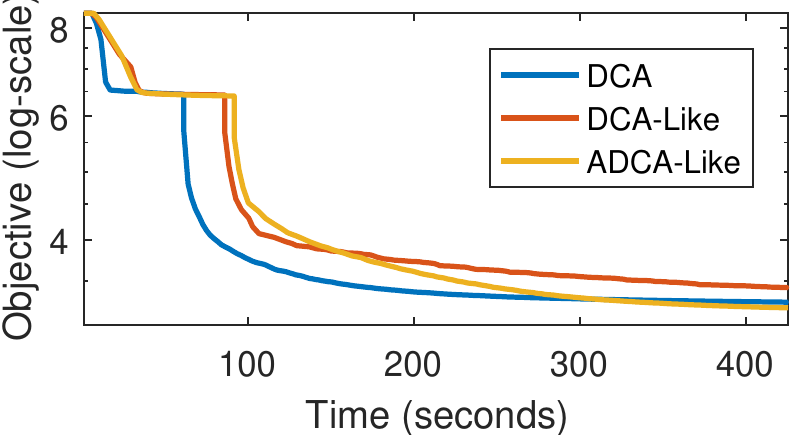} 
\\
(d) \textit{mnist} & (e) \textit{miniboone}  & (f) \textit{covertype} 
\\
\includegraphics[width=0.3\linewidth]{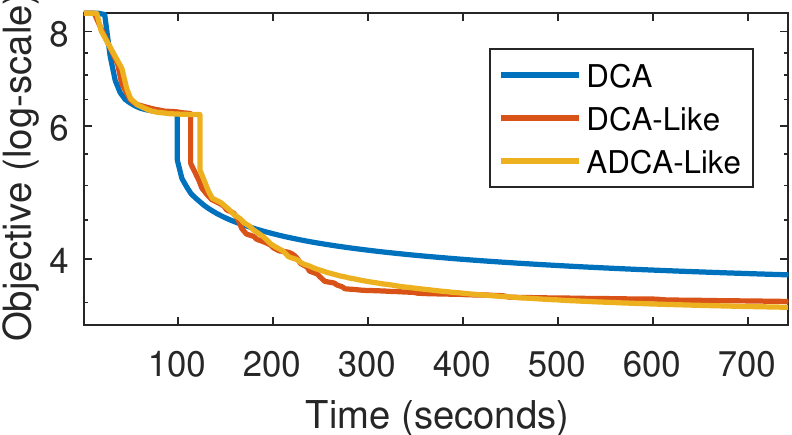} &
\includegraphics[width=0.3\linewidth]{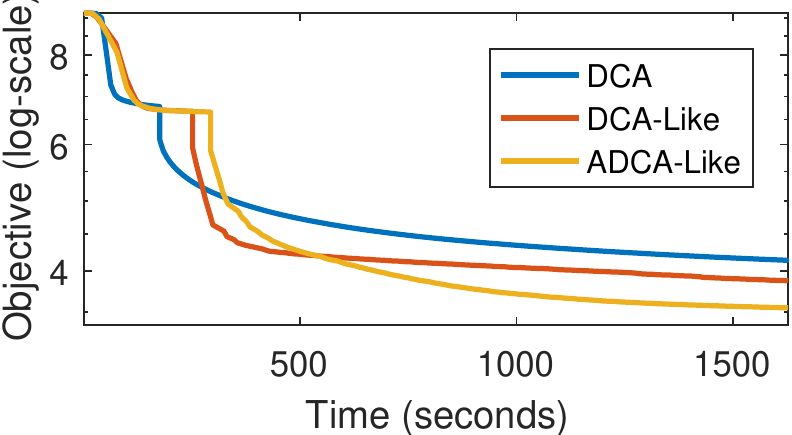} &
\includegraphics[width=0.3\linewidth]{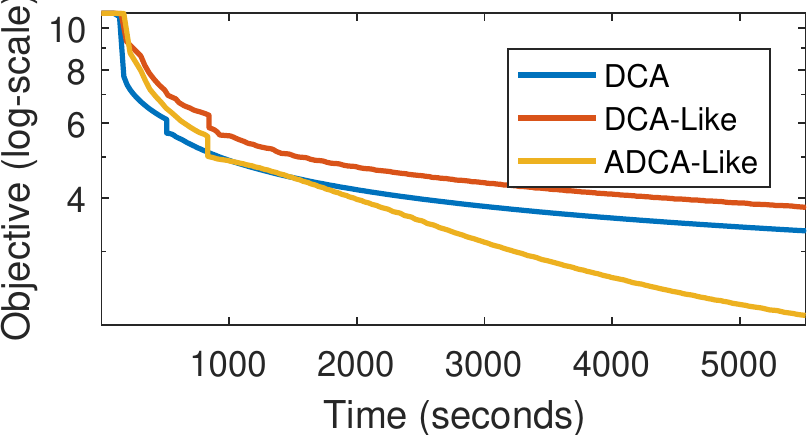}
\\
\end{tabular}
\captionof{figure}{Objective value versus running time (average of ten runs). 
}
\label{fig.plot_progress}
\end{table*}

\begin{table*}[tbh!]
\centering
\begin{tabular}{ccc}
\\
DCA  &
DCA-Like &
ADCA-Like
\\
\includegraphics[width=0.25\linewidth]{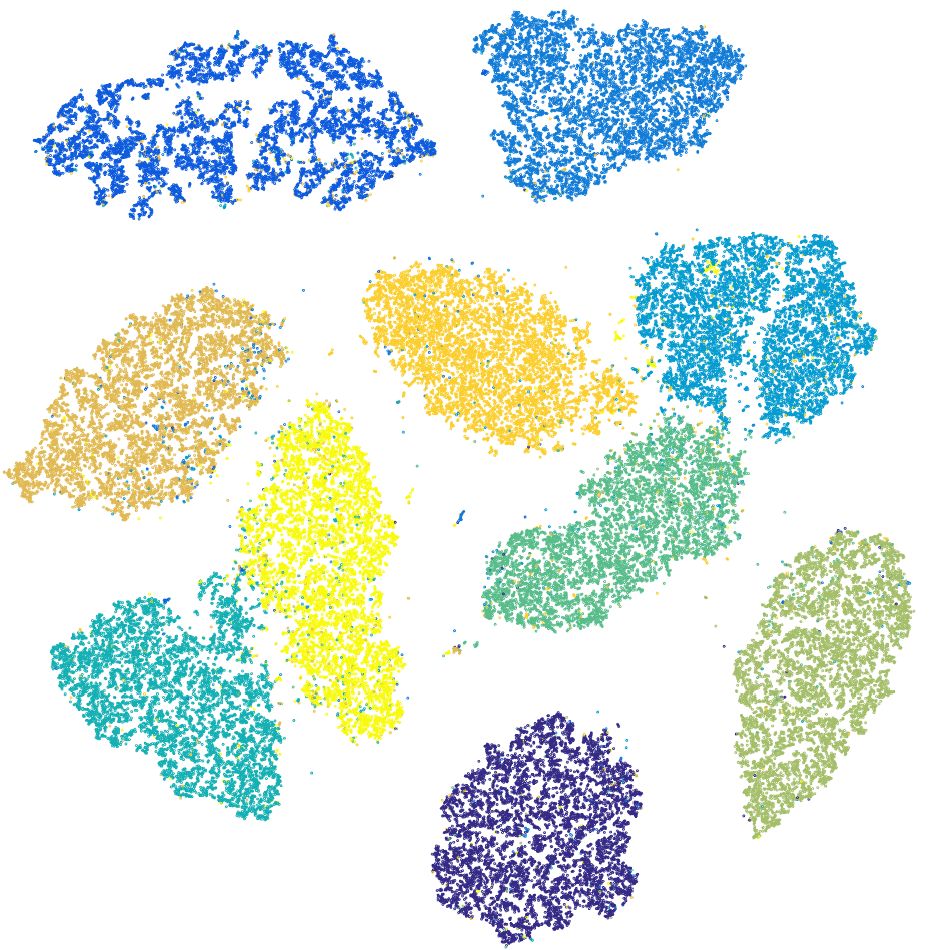} &
\includegraphics[width=0.25\linewidth]{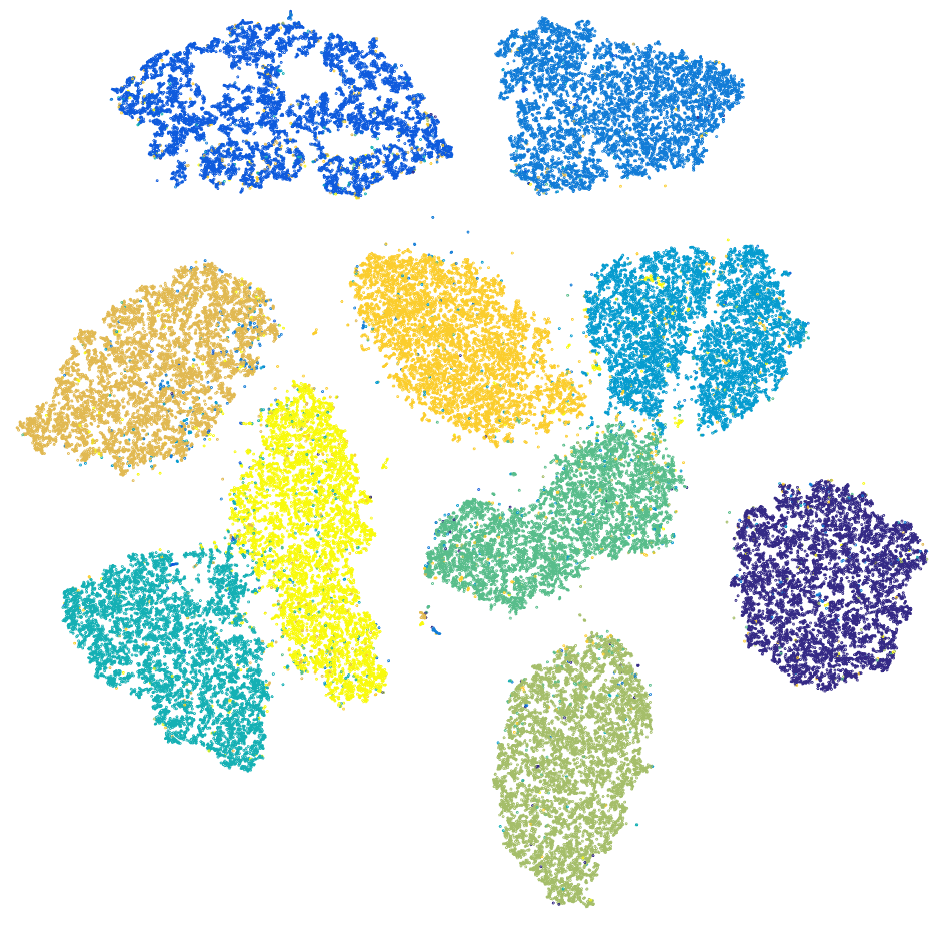} &
\includegraphics[width=0.25\linewidth]{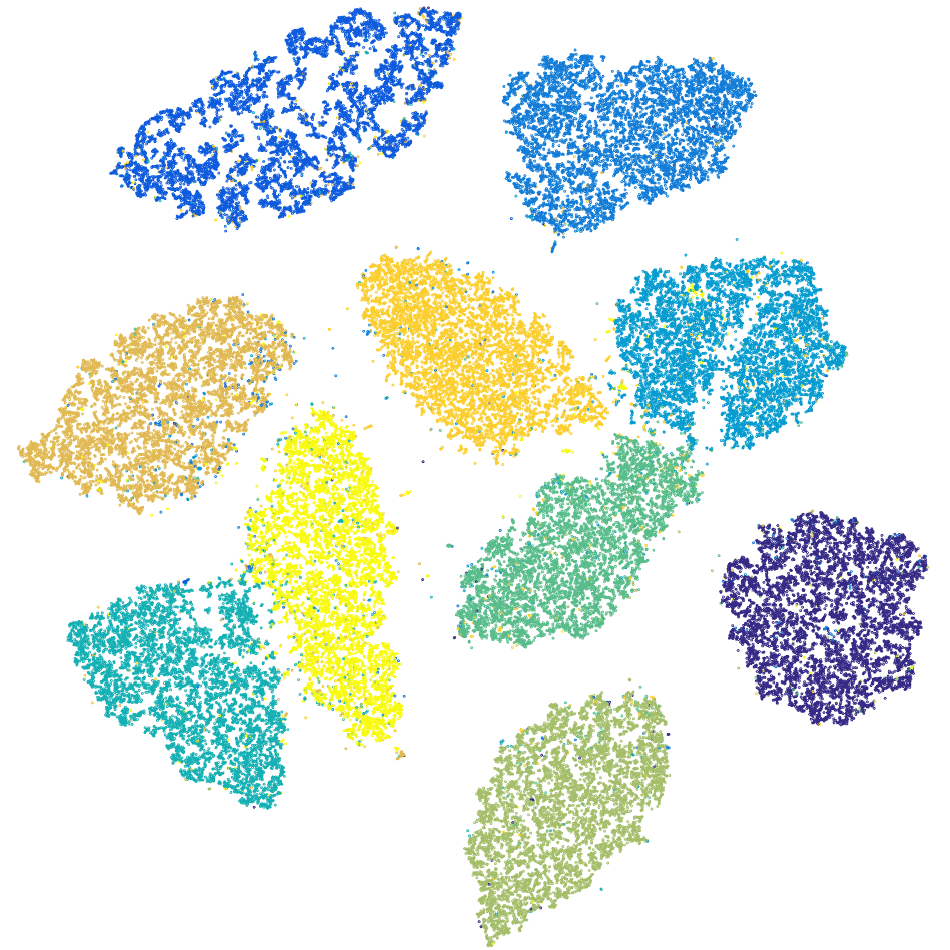}
\end{tabular}
\captionof{figure}{Visualization of embedding space on \textit{mnist} dataset. Colors represent classes of data (0-9).}
\label{fig.viz}
\end{table*}

Table \ref{tbl.exp_result} shows the average results of each algorithm after converging. DCA-Like is superior to DCA in all three criteria. In term of convergence speed, the number of iterations of DCA-Like is from 3 to 11.5 times less than DCA.  Consequently, the computing time of DCA-Like is improved by 1.9 to 9.2 times comparing to DCA. Furthermore, DCA-Like gives lower objective value than DCA in 4 out of 6 datasets (\textit{letters}, \textit{shuttle}, \textit{minibonne}, and \textit{covertype}), whereas the rest can be neglected.

ADCA-Like further improves the performance of DCA-Like. In term of number of iterations, ADCA-Like has the lowest by 1.8 to 5 times compared to DCA-Like. The gains in computing time are huge, as ADCA-Like is faster than DCA-Like (reps. DCA) from 1.5 to 5 times (reps. 6.5 to 13 times). Concerning objective, ADCA-Like performs the best in 5 out of 6 datasets among three algorithms. For only three cases (\textit{letters}, \textit{minibonne}, and \textit{mnist} dataset), DCA-Like performs as good as ADCA-Like, but for at least 1.5 times more time-consuming.

In Figure~\ref{fig.plot_progress}, we plot the value of objective function as time progress. Note that, we only plot the value until one of three algorithms stops. Surprisingly, DCA performs thoroughly at the beginning but then it is left behind; while both DCA-Like and ADCA-Like improves swiftly over time. It is noticeable that, in the plot of \textit{sensorless} and \textit{covertype}, DCA gives the better results than DCA-Like. Understandably, this Figure was captured at the end of ADCA-Like, when DCA-Like does not have enough time to surpass DCA but at the end DCA-Like gives better objective value than DCA (see Table~\ref{tbl.exp_result}).


Figure~\ref{fig.viz} visualizes \textit{mnist} dataset by all three algorithms.
This dataset consists of $70000$ gray-scale $28 \times 28$ images over $10$ classes of handwritten digits. \textit{mnist} can be considered as the benchmark dataset for SNE-based algorithms, since they are able to capture both local and global structure of this dataset, especially in 2D embedding space. 
As we can see, in the embedding space, all three algorithms managed to keep the structure of dataset on original space. Three images in Figure~\ref{fig.viz} are quite similar since the objective values of all algorithms in this case are fairly similar.

\section{Conclusions}\label{Sec:conclusion}
We have rigorously studied the the constrained sum of differentiable function and composite functions minimization problem. We reformulated the latter as a DC program and proposed two variants of DCA to solve the resulting problem.
In the first variant, we proposed a new technique to iteratively update the parameter $\mu$ and consequently the decomposition of objective function. We named the first variant as DCA-Like since the parameter $\mu$ is not large enough to ensure the successive decomposition of the objective function to be a DC decomposition. However, we proved that DCA-Like still enjoys the convergence properties of DCA. Furthermore, every limit point of the sequence generated by DCA-Like is a critical point. Considering the Kudyka-Lojasiewics assumption, we proved that each bounded sequence generated by DCA-Like globally converges to a critical point. The convergence rate under Kudyka-Lojasiewics assumption was also studied. In the second variant, ADCA-Like, we incorporate the Nesterov's acceleration technique into DCA-Like. We showed that ADCA-Like enjoys similar convergence properties and convergence rate of DCA-Like. To evaluate the performance of DCA-Like and ADCA-Like, we consider the t-distributed Stochastic Embedding (t-SNE) in data visualizing. DCA-Like and ADCA-Like applied to the t-SNE are inexpensive: the solution of the convex sub-problem can be explicitly computed. We showed that the Majorization-Minimization algorithm, the best state-of-the-art algorithm for t-SNE is nothing else but DCA-Like applied to t-SNE. Numerical experiments were carefully conducted on several benchmark datasets. The numerical results show that DCA-Like greatly improves the convergence speed of DCA. DCA-Like is up to 9.2 times faster than DCA while giving better objective value on 4 out of 6 datasets and comparable objective value on the 2 remaining datasets. ADCA-Likes improves further the convergence speed as well as the objective value of DCA-Like. The gain of computation time of ADCA-Like is up to 5 times smaller than DCA-Like. DCA-Like and ADCA-Like are undoubtedly improvements of DCA.

\bibliography{dca-tsne-icml18}
\bibliographystyle{icml2018}

%
%
%

\end{document}